\documentclass[reqno]{amsart}

\usepackage[T1]{fontenc}
\usepackage[utf8]{inputenc}

\usepackage{microtype}
\usepackage{lmodern}

\usepackage[unicode=true,pdfusetitle]{hyperref}
\let\C\undefined
\usepackage[abbrev]{amsrefs}

\numberwithin{equation}{section}
\newtheorem{theorem}{Theorem}[section]
\newtheorem{proposition}[theorem]{Proposition}

\theoremstyle{definition}
\newtheorem{definition}[theorem]{Definition}
\allowdisplaybreaks

\usepackage{enumitem}
\usepackage[abbrev]{amsrefs}

\usepackage{mathtools}
\usepackage{cleveref}

\newcommand{\Rset}{\mathbb{R}}

\newcommand{\Deriv}{\mathrm{D}}
\newcommand{\dext}{\mathrm{d}}
\newcommand{\Nset}{\mathbb{N}}
\newcommand{\Sset}{\mathbb{S}}

\newcommand{\defeq}{\coloneqq}

\newcommand{\dif}{\;\mathrm{d}}

\DeclarePairedDelimiter{\abs}{\lvert}{\rvert}
\DeclarePairedDelimiter{\norm}{\lVert}{\rVert}

\DeclarePairedDelimiter{\brk}{(}{)}

\DeclarePairedDelimiterX\dualprod[2]{\langle}{\rangle}{#1, #2}
\DeclarePairedDelimiterX\intvo[2]{(}{)}{#1, #2}
\DeclarePairedDelimiterX\intvc[2]{[}{]}{#1, #2}
\DeclarePairedDelimiterX\intvl[2]{(}{]}{#1, #2}
\DeclarePairedDelimiterX\intvr[2]{[}{)}{#1, #2}

\DeclarePairedDelimiterX\set[1]\{\}{%

#1
}
\DeclareMathOperator{\Lin}{Lin}

\crefname{subsection}{§}{§§}

\usepackage{constants}

\begin{document}
\title[Sobolev inequalities and cancelling operators]
{Endpoint Sobolev inequalities for vector fields and cancelling operators}

\author[Jean Van Schaftingen]{Jean Van Schaftingen}

\address{Jean Van Schaftingen\newline
Université catholique de Louvain (UCLouvain)\\
Institute de Recherche en Mathématique et Physique (IRMP)\\
Chemin du Cyclotron 2 bte L7.01.01\\
1348 Louvain-la-Neuve\\
Belgium}
\email{Jean.VanSchaftingen@UCLouvain.be}


\begin{abstract}
The injectively elliptic vector differential operators $A (\mathrm{D})$ from $V$ to $E$ on $\mathbb{R}^n$ such that the estimate
\[
 \Vert \mathrm{D}^\ell u\Vert_{L^{n/(n - (k - \ell))} (\mathbb{R}^n)}
 \le \Vert A (\mathrm{D}) u\Vert_{L^1 (\mathbb{R}^n)}
\]
holds can be characterized as the operators satisfying a cancellation condition
\[
 \bigcap_{\xi \in \mathbb{R}^n \setminus \{0\}} A (\xi)[V] = \{0\}\;.
\]
These estimates unify existing endpoint Sobolev inequalities for the gradient of scalar functions (Gagliardo and Nirenberg), the deformation operator (Korn–Sobolev inequality by M.J. Strauss) and the Hodge complex (Bourgain and Brezis).
Their  proof is based on the fact that $A (\mathrm{D}) u$ lies in the kernel of a cocancelling differential operator.
\end{abstract}

\maketitle 

\section{Sobolev inequalities for vector fields}

The classical \emph{Sobolev inequality} \citelist{\cite{Sobolev_1938}\cite{Gagliardo_1958}\cite{Nirenberg_1959}} states that given 
\(n, k \in \Nset \setminus \set{0}\) and \(p \in  [1, \frac{n}{k - \ell})\) there exists a constant \(\Cl{cst_Eeh9quaeQu2iyeeng1jaechu} \in \intvo{0}{\infty}\) such that each function \(u \in C^\infty_c (\Rset^n, \Rset)\) satisfies the inequality
\begin{equation}
\label{eq_AhR2ooxi0teixahheePoa5hu}
  \brk[\Big]{\int_{\Rset^n} \abs{\Deriv^{\ell} u}^\frac{np}{n - (k - \ell)p}}^{1 - \frac{(k - \ell)p}{n}}
  \le 
  \Cr{cst_Eeh9quaeQu2iyeeng1jaechu}
  \int_{\Rset^n} \abs{\Deriv^{k} u}^p\;.
\end{equation}
Given linear spaces \(V\) and \(E\) and a linear differential operator \(A (\Deriv)\) of order \(k \in \Nset \setminus \set{0}\) from \(V\) to \(E\) on \(\Rset^n\) defined for \(u \in C^\infty (\Rset^n, V)\) at each \(x \in \Rset^n\) by \(A (\Deriv) u(x) \defeq A(\Deriv u (x))\),  where \(A \in \Lin (\Lin_{\mathrm{sym}}^k (\Rset^n, V), E)\), or equivalently, by
\begin{equation*}
  A (\Deriv) u (x) = \sum_{\substack{\alpha \in \Nset^n\\ \abs{\alpha} = k}} A_\alpha [\partial^\alpha u (x)] = \sum_{\substack{\alpha \in \Nset^n\\ \abs{\alpha} = k}} \partial^\alpha (A_\alpha [u]) (x)\;,
\end{equation*}
with \(A_\alpha \in \Lin (V, E)\) for \(\alpha = (\alpha_1, \dotsc, \alpha_n) \in \Nset^n\) satisfying \(\abs{\alpha} \defeq \alpha_1 + \dotsb + \alpha_n = k\), 
the goal of the present work is to determine whether every vector field \(u \in C^\infty (\Rset^n, V)\) satisfies a vector Sobolev inequality 
\begin{equation}
  \label{eq_bah2acooqueBei4Ahph}
   \brk[\Big]{\int_{\Rset^n} \abs{\Deriv^{\ell} u}^\frac{np}{n - (k - \ell)p}}^{1 - \frac{(k - \ell)p}{n}}
  \le 
  \C
  \int_{\Rset^n} \abs{A (\Deriv)u}^p \;.
\end{equation}

When \(p \in \intvo{1}{\infty}\), the injective ellipticity is the key notion to have \eqref{eq_bah2acooqueBei4Ahph}.
\begin{definition}
   Given \(n \in \Nset \setminus \set{0}\), finite-dimensional vector spaces \(V\) and \(E\), a homogeneous constant coefficient differential operator \(A (\Deriv)\)
    of order \(k \in \Nset \setminus \set{0}\) from 
    \(V\) to \(E\) on \(\Rset^n\) 
    is \emph{injectively elliptic} whenever for every \(\xi \in \Rset^n\setminus \set{0}\), one has
    \(\ker A (\xi) = \set{0}\). 
\end{definition} 

If the operator \(A (\Deriv)\) is injectively elliptic, we can write in the Fourier domain for every \(\xi \in \Rset^n \setminus \set{0}\)
\begin{equation}
\label{eq_maiPe2pooshoo0vu8aecoy4a}
 \mathcal{F} (\Deriv^k u) (\xi)
 = \xi^{\otimes k} \otimes A (\xi)^\dagger[\mathcal{F} (A (\Deriv)u) (\xi)]\;,
\end{equation}
where the Fourier transform \(\mathcal{F} u : \Rset^n\to V + i V\) of a Schwartz test function \(u \in \mathcal{S} (\Rset^n, V)\) is defined for every \(\xi \in \Rset^n\) by the integral formula
\begin{equation*}
(\mathcal{F} u) (\xi) 
\defeq 
\int_{\Rset^n} e^{-2\pi i\,\xi\cdot x}  \, u (x) \dif x\;
\end{equation*}
and where for each \(\xi \in \Rset^n \setminus \set{0}\), \(A(\xi)^\dagger\) is the Moore--Penrose generalized inverse of \(A (\xi)\):
\begin{equation}
\label{eq_abihooJ4iateimeXie5ohy9f}
  A (\xi)^\dagger \defeq \brk*{A (\xi)^* A (\xi)}^{-1} A (\xi)^*\;,
\end{equation}
with \(A (\xi)^* \in \Lin (E, V)\) the adjoint of \(A(\xi)\).
Applying the Parseval identity when \(p = 2\) and a classical multiplier theorem (see for example \cite{Stein_1970}*{ch.\ IV th.\ 3}) when \(p \in \intvo{1}{\infty} \setminus \set{2}\),
 we obtain from \eqref{eq_maiPe2pooshoo0vu8aecoy4a} the estimate 
\begin{equation} 
\label{eq_chausaewagheich6eabizioS}
\int_{\Rset^n} \abs{\Deriv^{k} u}^p
  \le \C \int_{\Rset^n} \abs{A (\Deriv) u }^p \;.
\end{equation}
As a consequence of the inequalities \eqref{eq_chausaewagheich6eabizioS} and \eqref{eq_AhR2ooxi0teixahheePoa5hu}, we get that when \(p \in \intvo{1}{\frac{n}{k-\ell}}\), the estimate \eqref{eq_bah2acooqueBei4Ahph} holds when the operator \(A (\Deriv)\) is injectively elliptic.

The proof outlined above does not work at all in the endpoint \(p = 1\).
More dramatically, Ornstein 
\cite{Ornstein_1962} has proved that if \(B (\Deriv)\) 
is a homogeneous constant coefficient differential operator of order \(k\) from \(V\) to a linear space \(F\) on \(\Rset^n\) and if for every \(u \in C^\infty_c (\Rset^n, V)\) the estimate 
\begin{equation}
\label{eq_Eexiew6iequacahbahgaewei}
\int_{\Rset^n} \abs{B (\Deriv) u} \le \C \int_{\Rset^n} \abs{A (\Deriv)u}\;,
\end{equation}
holds, then one can write \(B (\Deriv) = L A (\Deriv)\), with \(L : E \to F\) being constant-coefficient linear mapping.

\section{Sobolev estimates for cancelling operators}

Continuing our investigation of the Sobolev-type inequality \eqref{eq_bah2acooqueBei4Ahph}, we will examine how badly \(A (\Deriv)u\) does not control \(u\) beyond Ornstein's non-estimate \eqref{eq_Eexiew6iequacahbahgaewei}.
In order to do this, we can try to have \(A (\Deriv) u\) as singular as possible, that is, close to a Dirac measure and thus to construct some \(u : \Rset^n \to V\) such that, for some fixed vector \(e \in E\),
\begin{equation}
\label{eq_ahfigh4ooPoh6eix7ohnii5Y}
 A (\Deriv) u = e \delta_0
\end{equation}
on \(\Rset^n\) in the sense of distributions. 
Taking the Fourier transform on both sides of \eqref{eq_ahfigh4ooPoh6eix7ohnii5Y}, 
we get for every \(\xi \in \Rset^n\)
\begin{equation}
\label{eq_Zahgh4Oob1ilaeyoh0shah4N}
 \brk{2 \pi i}^k A (\xi)[\mathcal{F} u (\xi)]= e\;.
\end{equation}
The equation \eqref{eq_Zahgh4Oob1ilaeyoh0shah4N} will have a solution for every \(\xi \in \Rset^n\) if and only if the operator \(A (\Deriv)\) does not satisfy the following cancellation condition, introduced by the author \cite{VanSchaftingen_2013}.

\begin{definition}
\label{definition_cancelling}
The homogeneous differential operator \(A (\Deriv)\) is \emph{cancelling} whenever 
\[
 \bigcap_{\xi \in \Rset^n \setminus \set{0}} A (\xi)[V] = \set{0}\;.
\]
\end{definition}

Thanks to the injective ellipticity of \(A (\Deriv)\), \(A (\xi)^\dagger\) can be defined by \eqref{eq_abihooJ4iateimeXie5ohy9f} and is homogeneous of degree \(-k\);
thanks to a construction \citelist{\cite{Bousquet_VanSchaftingen_2014}*{lem.\ 2.1}\cite{Raita_2019}} based on classical constructions in distribution theory and Fourier analysis \cite{Hormander_1990_I}*{th.\ 3.2.3, 3.2.4, 7.1.18, th.\ 7.1.16} (see also \cite{VanSchaftingen_2023}*{prop.\ 2} for a direct self-contained proof), 
one can construct a representation kernel \(G_A \in C^\infty (\Rset^n \setminus \set{0}, \Lin (E, V))\) such that for every \(\xi \in \Rset^n\),
\begin{equation*}
\mathcal{F} G_A (\xi)  = \brk*{2 \pi i}^{-k} A^\dagger (\xi) 
\end{equation*}
and such that 
for every \(x \in \Rset^n \setminus \set{0}\) and every \(t \in \Rset \setminus \set{0}\), 
\[
    G_A (t x) = t^{n - k} \brk[\big]{G_A (x) - \ln \abs{t} \, P_A (x)}\;,
\]
where the function \(P_A : \Rset^n \to \Lin (E, V)\) is a homogeneous polynomial of degree \(k - n\) when \(k - n \ge 0\) is even, and is \(0\) otherwise.

Thanks to \cref{definition_cancelling}, we can now state our main result characterizing endpoint Sobolev inequalities \cite{VanSchaftingen_2013}*{prop.\ 4.6 and  5.5}.

\begin{theorem}
\label{theorem_cancelling_necessary_Sobolev}
Let \(n \in \Nset\setminus\set{0}\), let \(V\) and \(E\) be finite-dimensional vector spaces 
and let \(A (\Deriv)\) be a homogeneous constant coefficient differential operator
of order \(k \in \Nset \setminus \set{0}\) from 
\(V\) to \(E\) on \(\Rset^n\).
If \(A (\Deriv)\) is injectively elliptic and if \(\ell \in \Nset \setminus \set{0}\) satisfies \(0 < k - \ell < n\), then there exists a constant \(\Cl{cst_yuquaichienoo4ReeZ5wus4a} \in \intvo{0}{\infty}\) such that for each \(u \in C^\infty_c (\Rset^n, V)\) 
\begin{equation}
\label{eq_eeghaib2eiviu0cib4Pei8ut}
  \brk[\Big]{\int_{\Rset^n} \abs{\Deriv^\ell u}^\frac{n}{n - (k - \ell)}}^{1 - \frac{k - \ell}{n}}
\le 
\Cr{cst_yuquaichienoo4ReeZ5wus4a}
\int_{\Rset^n} \abs{A (\Deriv)[u]}\;,
\end{equation}
if and only if \(A (\Deriv)\) is cancelling.
\end{theorem}

The necessity of the cancellation follows essentially by observing that if \(e \in \bigcap_{\xi \in \Rset^n \setminus \set{0}} A (\xi)[V] \setminus \set{0}\), then a suitable approximation of \(G_A [e]\) by smooth functions prevents \eqref{eq_eeghaib2eiviu0cib4Pei8ut} from holding.

A first consequence of \cref{theorem_cancelling_necessary_Sobolev}, is the endpoint 
Sobolev inequality of Bourgain and Brezis \citelist{\cite{Bourgain_Brezis_2004}\cite{Bourgain_Brezis_2007}*{cor.\ 17}} (see also \cite{Lanzani_Stein_2005}): given \(m \in \set{1, \dotsc, n - 1}\), 
for every \(u \in C^\infty_c (\Rset^n, \bigwedge^{m} \Rset^n)\),
\begin{equation}
\label{eq_shein0thah0fie7hoot6ic6E}
  \brk[\Big]{\int_{\Rset^n} \abs{u}^\frac{n}{n- 1}}^{1 - \frac{1}{n}}
  \le C \int_{\Rset^n} \abs{\dext u} + \abs{\dext^* u}
\end{equation}
holds if and only \(m \not \in \set{1, n - 1}\).
Here \(\dext u \in C^\infty_c (\Rset^n, \bigwedge^{m + 1} \Rset^n)\) and \(\dext^* u \in C^\infty_c (\Rset^n, \bigwedge^{m - 1} \Rset^n)\) denote respectively the exterior differential and codifferential of the differential form \(u\).

As a second consequence, we have Strauss's endpoint Korn--Sobolev inequality \cite{Strauss_1973}: for every \(u \in C^\infty_c (\Rset^n, \Rset^n)\), one has 
\begin{equation}
  \brk[\Big]{\int_{\Rset^n} \abs{u}^\frac{n}{n- 1}}^{1 - \frac{1}{n}}
  \le \C \int_{\Rset^n} \abs{\Deriv_{\mathrm{sym}} u}\;,
\end{equation}
where \(\Deriv_{\mathrm{sym}} u \defeq (\Deriv u + (\Deriv u)^*)/2\) is the \emph{symmetric derivative,} also known in elasticity as the \emph{deformation operator.}

\section{Duality estimates for cocancelling operators}

The proof of the sufficiency of the cancellation \eqref{theorem_cancelling_necessary_Sobolev} is based on the crucial fact that \(A (\Deriv) u\) on the right-hand side is not any function, but is constrained to satisfy some \emph{compatibility conditions}.

\begin{proposition}
  \label{proposition_compatibility_conditions}
Let \(n \in \Nset \setminus \set{0}\), let \(V\) and \(E\) be finite-dimensional vector spaces, and let \(A (\Deriv)\) be a homogeneous constant coefficient differential operator
of order \(k \in \Nset \setminus \set{0}\) from 
\(V\) to \(E\) on \(\Rset^n\).
If \(A (\Deriv)\) is injectively elliptic, then there exists a homogeneous constant coefficient differential operator \(L (\Deriv)\)
from 
\(E\) to \(E\) on \(\Rset^n\) such that for every \(\xi \in \Rset^n \setminus \set{0}\),
\begin{equation}
\label{eq_phicee9aad0ugh2thaeTueph}
  A (\xi)[V] = \ker L (\xi)\;.
\end{equation}
\end{proposition}

The proof of \cref{proposition_compatibility_conditions} is based on the definition of \(L (\Deriv)\) by requiring that for each \(\xi \in \Rset^n \setminus \set{0}\), 
\begin{equation}
 L (\xi) \defeq \det \bigl(A (\xi)^*A (\xi)\bigr)\,
 \brk[\big]{\operatorname{id}_E\, -\, A (\xi) (A (\xi)^* A (\xi))^{-1} A (\xi)^*}\;.
\end{equation}

\Cref{proposition_compatibility_conditions} can be seen as a generalization of the \emph{symmetry of second-order derivatives}
\begin{equation}
\label{eq_eiy4ooXahk0Ood6atha3ohC4}
  \partial_j (\partial_i u) = \partial_i (\partial_j u)
\end{equation}
and of the \emph{Saint-Venant compatibility conditions} for the symmetric derivative
\begin{equation}
\label{eq_sahPhai2aeg2eo4moMee6woo}
 \partial_{k\ell} (\partial_i u^j +\partial_j u^i)
 +\partial_{ij} (\partial_k u^\ell +\partial_\ell u^k)
 = \partial_{kj} (\partial_i u^\ell +\partial_\ell u^i)
 +\partial_{i\ell} (\partial_k u^j +\partial_j u^k)\;,
\end{equation}
although the construction of \cref{proposition_compatibility_conditions} gives a more complicated operator than what appears in \eqref{eq_eiy4ooXahk0Ood6atha3ohC4} and \eqref{eq_sahPhai2aeg2eo4moMee6woo}.

The definition of cancelling operator (\cref{definition_cancelling}) and the construction of compatibility conditions (\cref{proposition_compatibility_conditions}) suggest the definition of cocancelling operator.

\begin{definition}
  \label{definition_cocancelling}
  Let \(n \in \Nset \setminus \set{0}\) and let \(E\) and \(F\) be finite-dimensional vector spaces.
A homogeneous constant coefficient differential operator \(L(\Deriv)\) from \(E\) to \(F\) on \(\Rset^n\) is \emph{cocancelling} whenever 
\begin{equation*}
\bigcap_{\xi \in \Rset^n \setminus \set{0}} \ker L (\xi) = \set{0}\;. 
\end{equation*}
\end{definition}

The cocancellation condition characterizes the operators for which there is a duality estimate with critical Sobolev spaces \citelist{\cite{VanSchaftingen_2008}\cite{VanSchaftingen_2013}} (see also previous results \citelist{\cite{Bourgain_Brezis_Mironescu_2004}\cite{VanSchaftingen_2004_circ}\cite{Bourgain_Brezis_2004}\cite{VanSchaftingen_2004_div}\cite{Bourgain_Brezis_2007}\cite{VanSchaftingen_2004_ARB}}).

\begin{theorem}
\label{theorem_cocancelling}
  Let \(n \in \Nset \setminus \set{0, 1}\), let \(V\) and \(E\) be finite-dimensional vector spaces, let \(L (\Deriv)\) be a homogeneous constant coefficient differential operator from \(E\) to \(F\) on \(\Rset^n\) and let \(\ell \in \set{1, \dotsc, n - 1}\). 
  There exists a constant \(\Cl{cst_yoo5aek4Loo7eiyuZeiv5Ier} \in \intvo{0}{\infty}\) such that
for every \(f \in L^1 (\Rset^n, E)\) that satisfies \(L (\Deriv) f = 0\) in \(\Rset^n\) in the sense of distributions and every \(\varphi \in C^\infty_c (\Rset^n, E)\) one has 
\begin{equation}
  \label{eq_koo4lui9Eipo7XeenaeghaoR}
  \abs[\Big]{\int_{\Rset^n} \dualprod{f}{\varphi}}
  \le 
  \Cr{cst_yoo5aek4Loo7eiyuZeiv5Ier} \int_{\Rset^n} \abs{f}\; \brk[\Big]{\int_{\Rset^n} \abs{\Deriv^\ell \varphi}^\frac{n}{\ell}}^\frac{\ell}{n}
\end{equation}
if and only if the operator \(L (\Deriv)\) is cocancelling.
\end{theorem}

\Cref{theorem_cocancelling} states somehow that with regards to integration against vector fields that are in the kernel functions in the homogeneous Sobolev space \(\dot{W}^{\ell, n/\ell} (\Rset^n, E)\) behave as if they were bounded --- which is well-known not to be the case.

The necessity of the cancellation can be seen by noting that if \(L (\Deriv)\) was not cancelling, then one would have \(L(\Deriv) (\delta_0 e)= 0\) for some \(e \in E \setminus \set{0}\); approximating the measure \(\delta_0 e\) by smooth functions one would deduce from \eqref{eq_koo4lui9Eipo7XeenaeghaoR} that the Sobolev space \(\dot{W}^{\ell, n/\ell} (\Rset^n, \Rset)\) would be continuously embedded in \(L^\infty (\Rset^n, \Rset)\), which is not the case when \(\ell \in \set{1, \dotsc, n - 1}\).

Assuming that \cref{theorem_cocancelling}, \cref{theorem_cancelling_necessary_Sobolev} can be proved as follows. Noting that the operator \(L(\Deriv)\) given by \cref{theorem_cocancelling} is cocancelling, one gets by \cref{theorem_cocancelling} that \(\norm{A (D) u}_{W^{-\ell, n/(n- \ell)}} \le\C \norm{A (D) u}_{L^1}\) so that a classical multiplier theorem brings the conclusion.

Bourgain and Brezis's original proof \cite{Bourgain_Brezis_2007} of estimates of the type of \cref{theorem_cocancelling} was based on an approximation property for critical Sobolev functions through a Littlewood--Paley decomposition, generalizing a similar results in the study of the divergence equation \citelist{\cite{Bourgain_Brezis_2002}\cite{Bourgain_Brezis_2003}}; the advantage of their proof compared to the one presented below is that it provides much stronger estimates of the form 
\begin{equation}
  \abs[\Big]{\int_{\Rset^n} f \cdot \varphi}
  \le \C \norm{f}_{L^1 (\Rset^n) + \dot{W}^{-1,n/(n-1)}(\Rset^n)}
  \brk[\bigg]{\int_{\Rset^n} \abs{\Deriv \varphi}^n}^\frac{1}{n}\;.
\end{equation}

Let us now explain how \cref{theorem_cocancelling} can be proved in the case where \(L(\Deriv)\) is the divergence operator, following \cite{VanSchaftingen_2004_div}. (The reader is referred to \citelist{\cite{VanSchaftingen_2008}\cite{VanSchaftingen_2013}\cite{VanSchaftingen_2023}\cite{VanSchaftingen_2023_Paseky}} for the general case.)
Without loss of generality, we are going to estimate the integral
\begin{equation}
\label{eq_AehaeZ4dieh7ichax4ieyaex}
  \int_{\Rset^n} f \cdot e_n\, \phi\;,
\end{equation}
for \(\phi \in C^\infty_c (\Rset^n, \Rset)\), where \(e_n\) is the \(n^\text{th}\) vector in the canonical basis of \(\Rset^n\).
By Fubini's theorem we have
\begin{equation}
\label{eq_Ohseequaz4Faing1ingeeshi}
  \int_{\Rset^n} f \cdot e_n\, \phi
  = \int_{\Rset} \brk[\Big]{\int_{\Rset^{n - 1}} f  (\cdot, x_n) \cdot e_n \,\phi (\cdot, x_n)} \dif x_n\;.
\end{equation}
We are now going to estimate the inner integral on the right-hand side of \eqref{eq_Ohseequaz4Faing1ingeeshi}.
First, we immediately have for each \(\psi \in C^\infty (\Rset^{n - 1}, \Rset)\),  
\begin{equation}
\label{eq_tie1Lefemieyum2zasoeWaip}
 \abs[\Big]{\int_{\Rset^{n - 1}} f  (\cdot, x_n) \cdot e_n\, \psi} 
 \le\norm{\psi}_{L^\infty (\Rset^{n - 1})} \int_{\Rset^{n - 1}} \abs{f(\cdot, x_n)} \;.
\end{equation}
On the other hand, by the Gauss--Ostrogradsky divergence theorem,  we also have
\begin{equation}
\label{eq_eiqu4daeFuDeiv0ZaiPuijai}
 \abs[\Big]{\int_{\Rset^{n - 1}} f  (\cdot, x_n) \cdot e_n\, \psi}
 = \abs[\Big]{\int_{\Rset^{n - 1} \times \intvo{x_n}{\infty}} \operatorname{div} (f \Psi)}
 \le \norm{\Deriv \psi}_{L^\infty (\Rset^{n - 1})} \int_{\Rset^n} \abs{f} \;. 
\end{equation}
since 
\(\operatorname{div} f = 0\).
Interpolating between the estimates \eqref{eq_tie1Lefemieyum2zasoeWaip} and \eqref{eq_eiqu4daeFuDeiv0ZaiPuijai} and applying the Morrey--Sobolev embedding on \(\Rset^{n - 1}\), we get
\begin{equation}
\label{eq_eipiu2waeTh8ooboophio5te}
\begin{split}
  \abs[\Big]{\int_{\Rset^{n - 1}} f  (\cdot, x_n) \cdot e_n\, \psi}
  &\le \C \brk[\Big]{\int_{\Rset^{n - 1}} \abs{f(\cdot, x_n)}}^{1 - 1/n} \brk[\Big]{\int_{\Rset^n} \abs{f}}^{1/n} \abs{\psi}_{C^{0, 1/n} (\Rset^{n - 1})}\\
  &\le \Cl{cst_Ahgh1GewophaFoh1laiCuiph} \brk[\Big]{\int_{\Rset^{n - 1}} \abs{f(\cdot, x_n)}}^{1 - \frac{1}{n}} \brk[\Big]{\int_{\Rset^n} \abs{f}}^\frac{1}{n} \brk[\Big]{\int_{\Rset^{n - 1}} \abs{\Deriv \psi}^n}^\frac{1}{n}\;.
\end{split}
\end{equation}
Combining \eqref{eq_Ohseequaz4Faing1ingeeshi} and \eqref{eq_eipiu2waeTh8ooboophio5te}, we deduce in view of Hölder's inequality that 
\begin{equation}
\label{eq_ieBaecheinohth1ohza1OhQu}
\begin{split}
 \abs[\Big]{\int_{\Rset^{n}} f \cdot e_n\, \phi}
 &\le \Cr{cst_Ahgh1GewophaFoh1laiCuiph}
\brk[\Big]{\int_{\Rset^n} \abs{f}}^\frac{1}{n}
 \int_{\Rset} \brk[\Big]{\int_{\Rset^{n - 1}} \abs{f(\cdot, x_n)}}^{1 - \frac{1}{n}} \brk[\Big]{\int_{\Rset^{n - 1}} \abs{\Deriv \phi (\cdot, x_n)}^n}^\frac{1}{n} \dif x_n\\
 &\le \Cr{cst_Ahgh1GewophaFoh1laiCuiph}
 \int_{\Rset^n} \abs{f} \, \brk[\Big]{\int_{\Rset^{n}} \abs{\Deriv \phi}^n}^\frac{1}{n} \;.
\end{split}
\end{equation}
Therefore for any vector \(\nu \in \Rset^n\), it follows from \eqref{eq_ieBaecheinohth1ohza1OhQu} that we have proved
\begin{equation}
\label{eq_woNgoo3heiD5Oogh8eneeguo}
 \abs[\Big]{\int_{\Rset^{n}} f \cdot \nu\, \phi}
 \le \Cr{cst_Ahgh1GewophaFoh1laiCuiph}\abs{\nu}
 \brk[\Big]{\int_{\Rset^n} \abs{f}} 
 \brk[\Big]{\int_{\Rset^{n}} \abs{\Deriv \phi}^n}^\frac{1}{n} \;.
\end{equation}
Decomposing \(\varphi = \sum_{j = 1}^n e_i \phi_i\) with \(\phi_i \defeq e_i \cdot \varphi\), we finally obtain \eqref{eq_koo4lui9Eipo7XeenaeghaoR} from \eqref{eq_woNgoo3heiD5Oogh8eneeguo}.

\Cref{theorem_cocancelling} when \(L (\Deriv)\) is the divergence 
is equivalent, thanks to Smirnov’s result on the approximation of divergence-free measures \cite{Smirnov_1993}, to the Bourgain, Brezis and Mironescu’s estimate on circulation integrals \cite{Bourgain_Brezis_Mironescu_2004} (see also \cite{VanSchaftingen_2004_circ}): if \(\Gamma \subset \Rset^n\) is a closed curve with tangent vector \(t\) and length \(\abs{\Gamma}\), then for every vector field \(\varphi \in C^\infty_c (\Rset^n, \Rset^n)\), one has
\begin{equation}
\label{eq_aekaemoovuphaich7ohGhie3}
 \abs[\Big]{\int_{\Gamma} \dualprod{\varphi}{t}}
 \le \C \abs{\Gamma} \brk[\Big]{\int_{\Rset^n} \abs{\Deriv \varphi}^n}^\frac{1}{n}\;;
\end{equation}
the geometric flavour of \eqref{eq_aekaemoovuphaich7ohGhie3} raises several natural open questions on sharp constants in higher dimensions \(n \ge 3\) \cite{Brezis_VanSchaftingen_2008}.

\section{Further results}
The cancellation condition can also be proved to be a necessary and sufficient condition for other estimates.

In the scale of fractional Sobolev spaces, for any injectively elliptic operator \(A (\Deriv)\) and assuming that \(k, \ell \in \Nset\), \(p \in \intvo{1}{\infty}\) and \(\sigma \in \intvo{0}{1}\) satisfy \(k - n = \ell + \sigma - \frac{n}{p}\), the estimate 
\begin{equation}
\label{eq_thoh2gohz4phahquahvohSuy}
\brk[\Big]{\int_{\Rset^n} \int_{\Rset^n}  \frac{\abs{\Deriv^{\ell} u (y) - \Deriv^{\ell} u (x)}^p} {\abs{y -x}^{n + \sigma p}} \dif y \dif x}^\frac{1}{p} 
\le 
\C
\int_{\Rset^n} \abs{A (\Deriv)[u]}\;,
\end{equation}
holds for every \(u \in C^\infty_c (\Rset^n, V)\) if and only if 
if and only if the operator \(A (\Deriv)\) is cancelling \cite{VanSchaftingen_2013} (see also \cite{VanSchaftingen_2023}).

Similarly, for any injectively elliptic operator \(A (\Deriv)\), the Hardy inequality 
\begin{equation}
\label{eq_ahHaichooxah6pha5thah6ie}
\int_{\Rset^n} \frac{\abs{\Deriv^\ell u(x)}}{\abs{x}^{k - \ell}} \dif x 
\le 
\C
\int_{\Rset^n} \abs*{A (\Deriv)[u]}\;,
\end{equation}
holds for every \(u \in C^\infty_c (\Rset^n, V)\) if and only if 
if and only if \(A (\Deriv)\) is cancelling \cite{VanSchaftingen_2013} (see also \cite{VanSchaftingen_2023} for the proof); this results originates in Maz\cprime{}ya’s work \cite{Mazya_2010} (see also \cite{Bousquet_Mironescu_2011}).

Finally, Raiță \cite{Raita_2019} has proved that for any injectively elliptic operator \(A (\Deriv)\) the uniform estimate 
\begin{equation}
\label{eq_cequae1eogheG4tohHohB9ae}
\norm{\Deriv^{n - k} u}_{L^\infty (\Rset^n)}
\le 
\C
\int_{\Rset^n} \abs{A (\Deriv)[u]}\;,
\end{equation}
is equivalent to the \emph{weak cancellation} property that for every \(e \in E\)
\begin{equation}
  \label{eq_aiboh2iesaeMisee2oyahah4}
  \int_{\Sset^{n - 1}} \xi^{\otimes k - n} A^{-1}(\xi) [e] \dif \xi = 0\;
\end{equation}
(see also \cite{VanSchaftingen_2023}*{\S 5.4}).

Endpoint estimates similar to \cref{theorem_cancelling_necessary_Sobolev,theorem_cocancelling} can also be obtained on stratified homogeneous groups \citelist{\cite{VanSchaftingen_Yung_2022}\cite{Chanillo_VanSchaftingen_2009}}, on the hyperbolic plane \cite{Chanillo_VanSchaftingen_Yung_2017_Variations} and on symmetric spaces of noncompact type \cite{Chanillo_VanSchaftingen_Yung_2017_Symmetric}.

For a more detailed exposition on endpoint Sobolev inequalities and cancelling operators, we refer the reader to the quite formal lecture notes \cite{VanSchaftingen_2023}, their somehow more informal counterpart \cite{VanSchaftingen_2023_Paseky} and to the survey article \cite{VanSchaftingen_2014}.

\end{document}